\documentclass[11pt]{article}
\usepackage{latexsym}
\usepackage{amssymb}

\newtheorem{teo}{THEOREM}[section]
\newtheorem{prop}[teo]{PROPOSITION}
\newtheorem{lem}[teo] {LEMMA}
\newtheorem{ejem}[teo]{EXAMPLE}
\newtheorem{obser}[teo]{Remark}
\newtheorem{defi}[teo]{DEFINITION}
\newtheorem{coro}[teo]{COROLLARY}
\newenvironment{proof}{\noindent \sf Proof. \sf}

\def\ker{\mathop{\rm Ker}\nolimits}

\def\Mod{\mathop{\rm Mod}\nolimits}

\def\cc{\mathcal{C}}
\def\dd{\mathcal{D}}
\def\bb{\mathcal{B}}
\def\M{\mathcal{M}}
\def\N{\mathcal{N}}

\begin{document}
\sf

\title{Galois coverings, Morita equivalence and smash extensions of categories over a field\thanks{This
work has been supported by the projects CONICET-CNRS:"METODOS
HOMOLOGICOS EN REPRESENTACIONES Y ALGEBRA DE HOPF", PICS 1514,
PICT 08280 (ANPCyT), UBACYTX169 and PIP-CONICET 02265. The second
author is a research member of CONICET (Argentina) and a Regular
Associate of ICTP Associate Scheme.}}
\author{Claude Cibils and Andrea Solotar}

\date{}

\maketitle

\begin{abstract}
We consider categories over a field $k$ in order to prove that
smash extensions and Galois coverings with respect to a finite
group coincide up to Morita equivalence of $k$-categories. For
this purpose we describe processes providing Morita equivalences
called contraction and expansion. We prove that composition of
these processes provides any Morita equivalence, a result which is
related with
 the karoubianisation (or idempotent completion) and additivisation of a $k$-category.
\end{abstract}



\font\fivrm=cmr5 \relax

\section {\sf Introduction }

In this paper we consider categories $\cc$ over a field $k$, which
means that the objects $\cc_0$ are a set, each morphism set
${}_y\cc_x$ from an object $x$ to an object $y$ is a $k$-vector
space and the composition of maps of $\cc$ is $k$-bilinear. In
particular each endomorphism set ${}_x\cc_x$ is an associative
$k$-algebra. Such categories are called $k$-categories, they have
been considered extensively and are considered as algebras with
several objects, see \cite{mi1,mi2}.

This work has a two-fold main purpose. In one direction we show
that there is a coincidence up to Morita equivalence between
 Galois coverings of
$k$-categories and smash extensions for a finite group. More
precisely we associate to each Galois covering of a $k$-category
with finite group $G$ a smash extension with the same group,
having the property that the categories involved are Morita
equivalent to the starting ones.  In particular from a full and
dense functor we obtain a faithful one. Conversely, a smash
extension of categories gives rise to a Galois covering, with
categories actually equivalent to the original ones. Consequently
both procedures are mutual inverses up to Morita equivalence.

In the other direction we study the Morita equivalence of
$k$-categories that we need for the results stated above. We
consider
 modules over a $k$-category $\cc$, that is $k$-functors from
$\cc$ to the category of $k$-vector spaces \emph{i.e.} collections
of vector spaces attached to the objects with "actions" of
morphisms transforming vectors at the source of the morphism to
vectors at the target. Notice that if $\cc$ is a finite object set
$k$-category it is well known and easy to prove that modules over
$\cc$ coincide with usual modules over the "matrix algebra"
$a(\cc)=\oplus_{x,y\in\cc}{}_y\cc_x$.

We introduce in this paper a general framework for Morita theory
for $k$-categories. More precisely we establish processes which
provide categories  Morita equivalent to a starting one. We prove
in the Appendix that up to equivalence of categories any Morita
equivalence of $k$-categories is a composition of contractions and
expansions of a given $k$-category, where contraction and
expansion are processes generalizing a construction considered in
\cite{cire}. More precisely, given a partition $E$ of the set of
objects of a $k$-category $\cc$ by means of finite sets, the
\emph{contracted category} $\cc_E$ along $E$ has set of objects
the sets of the partition while morphisms are provided by the
direct sum of all the morphism spaces involved between two sets of
the partition. The reverse construction is called
\emph{expansion}. Another process is related to the classical
Morita theory for algebras, that is for each vertex we provide an
endomorphism algebra Morita equivalent to the given one together
with a corresponding Morita context, which enables us to modify
the morphisms of the original category. In particular the matrix
category of a given category is obtained in this way. A discussion
of this processes in relation with karoubianisation and
additivisation (see for instance \cite{basc,ti}) is also presented
in the Appendix. We thank Alain Brugui\`{e}res and Mariano Suarez
Alvarez for useful conversations concerning this point.

Usually smash extensions are considered for algebras, see for
instance \cite{mo}. We begin by extending this construction to
$k$-categories, namely given a Hopf algebra $H$ we consider a Hopf
module structure on a $k$-category $\cc$ which is provided by an
$H$-module structure on each morphism space such that the
composition maps of $\cc$ are $H$-module maps - in particular  the
endomorphism algebra of each object is required to be an usual
$H$-module algebra. Given a Hopf module $k$-category $\cc$ we
 define the \emph{smash category} ${\cc}\#H$ in a
coherent way with the algebra case.

We need this extension of the usual algebra setting to the
categorical one in order to relate smash extensions to Galois
coverings of $k$-categories as considered for instance in
\cite{boga,cire,drov}.

Notice that we can consider, as in the algebra case, a smash
extension of a category as a Hopf Galois extension with the normal
basis property and with trivial map $\sigma$, see \cite[p.
101]{mo} and also \cite{beto,kr,rasa}. It would be interesting to
relate non trivial maps $\sigma$ to an extended class of coverings
of categories accordingly, we will not initiate this study in the
present paper.

We define a \emph{smash extension} of an $H$-module category $\cc$
to be the natural functor from $\cc$ to ${\cc}\#H$. An expected
compatibility result holds, namely if the number of objects of
$\cc$ is finite, the corresponding matrix algebra $a(\cc)$ has an
usual smash extension provided by $a(\cc\#H)$. The later algebra
can indeed be considered since the category ${\cc}\#H$ has also a
finite number of objects, namely the set of objects of $\cc$.
Moreover, we have that $a(\cc)\#H = a\left(\cc\#H\right)$.

We consider also Galois coverings of $k$-categories given by a
group $G$, that is a $k$-category with a free $G$-action and the
projection functor to the corresponding quotient category. More
precisely, by definition a $G$-$k$-category $\cc$ has a set action
of $G$ on the set of objects, and has linear maps ${}_y\cc_x
\to{}_{sy}\cc_{sx}$ for each element $s$ of $G$ and each couple of
objects $x$ and $y$, verifying the usual axioms that we recall in
the text. In other words we have a group morphism from $G$ to the
autofunctors of $\cc$. In case $\cc$ is object-finite, we infer a
usual action of $G$ by automorphisms of the algebra $a(\cc)$. A
$G$-$k$-category is called \emph{free} in case the set action on
the objects is free, namely $sx=x$ implies $s=1$. The quotient
category is well defined only in this case and we recall its
construction, see \cite{boga,ga,drov,cire,cima}.

The group algebra $kG$ is a Hopf algebra, hence we can consider
$kG$-module categories. Notice that $G$-$k$-categories form a
wider class than $kG$-module categories. In fact $kG$-module
categories are $G$-$k$-categories which have trivial action of $G$
on the set of objects.

First we establish a comparison between two constructions obtained
when starting with a graded category $\cc$ over a finite group
$G$. From one side the smash product category $\cc\# k^G$ is
defined in the present paper, and from the other side a smash
 product category $\cc\# G$ has been considered in
\cite{cima}, actually the later is the Galois covering of $\cc$
corresponding to the grading. We show that $\cc\# k^G$ and $\cc\#
G$ are not equivalent but Morita equivalent categories.

We note that starting with a Galois covering $\cc$ of a category
$\bb$, the covering category $\cc$ is $\bb\# G$ (see \cite{cima}
and the grading of $\bb$ introduced there, first considered by E.
Green in \cite{gr} for presented $k$-categories by a quiver with
relations). Unfortunately $\bb\# G$ has no natural $kG$-module
category structure. However $\bb\# G$ and $\bb\# k^G$ are Morita
equivalent and we perform the substitution. The later is a
$kG$-module category using the left $kG$-module structure of $k^G$
provided by $t\delta_s=\delta_{st^{-1}}$. In this way we associate
to the starting Galois covering the smash extension $(\bb\# k^G)\#
kG$ of $\bb\# k^G$.

The important point is that the later is Morita equivalent to
$\cc$ while $(\bb\# k^G)\# kG$ is isomorphic to a  matrix category
that we introduce, which in turn is Morita equivalent to $\bb$.
Notice that this result is a categorical version of the Cohen
Montgomery duality Theorem, see \cite{como}. Hence we associate to
the starting Galois covering $\cc\to\bb$  a smash extension with
the same group and where the categories are replaced by Morita
equivalent ones.

Second we focus to the reverse procedure, namely given a smash
extension of categories with finite group $G$ -- that is a
$kG$-module category $\bb$ and the inclusion $\bb\to\cc=\bb\#kG$
-- we intend to associate a Galois covering to this data. For this
purpose we consider the \emph{inflated category} $I_F\bb$ of a
category $\bb$ along a sequence $F=\{F_x\}$ of sets associated to
the vertices of the original category : each object $x$ of $\bb_0$
provides $\mid F_x\mid$ new objects while the set of morphisms
from $(x,i)$ to $(y,j)$ is precisely the vector space ${}_y\bb_x$
with the obvious composition. For a finite group $G$ the inflated
category of a $kG$-module category -- using the constant sequence
of sets $G$ -- has a natural structure of a free $G$-$k$-category.
The inflated category $I_{G}\bb$ is Morita equivalent to the
matrix category $M_{\mid G\mid}(\bb)$ by contraction and in turn
the later is Morita equivalent to $\bb$.

Moreover the categorical quotient of $I_{G}\cc$ exists and in this
way we obtain a Galois covering having the required properties
with respect to the starting smash extension.

\section{\sf Hopf module categories}

In this section we introduce the smash product of a category with
a Hopf algebra and we specify this construction in case the Hopf
algebra is the function algebra of a finite group $G$. We will
obtain that the later is Morita equivalent to the smash product
category defined in \cite{cima}.

We recall (see for instance \cite{mo}) that for a Hopf algebra $H$
over $k$, an $H$-module algebra $A$ is a $k$-algebra which is
simultaneously an $H$-module in such a way that the product map of
$A$ is a morphism of $H$-modules, where $A \otimes A$ is
considered as an $H$-module through the comultiplication of $H$.
Moreover we require that $h1_A = \epsilon(h)1_A$ for every $h\in
H$.

We provide an analogous definition for a $k$-category $\cc$
instead of an algebra.
\begin{defi}
A $k$-category $\cc$ is an $H$-\emph{module category} if each
morphism space is an $H$-module, each endomorphism algebra is an
$H$-module algebra and composition maps are morphisms of
$H$-modules, where as before the tensor product of $H$-modules is
considered as an $H$-module via the comultiplication of $H$.
\end{defi}
Notice that analogously we may consider the structure of an
$H$-comodule category. In case $H$ is a finite dimensional Hopf
algebra, we recall from \cite{mo} that there is a bijective vector
space preserving correspondence between right $H$-modules and left
$H^*$-comodules.

\begin{obser}
Given a finite $k$-category $\cc$, let $a(\cc)$ be the $k$-algebra
obtained as the direct sum of all $k$-module morphisms of $\cc$
equipped with the usual matrix product combined with the
composition of $\cc$. In case $\cc$ is an $H$-module category
$a(\cc)$ is an $H$-module algebra.
\end{obser}

Let $\cc$ be an $H$-module category. We define the  $k$-category
$\cc \# H$ as follows. The objets remain the same, while given two
objects $x$ and $y$ we put $_y(\cc \# H)_x = {}_y \cc_x \otimes_k
H$. The composition map for morphisms
$$_z (\cc \# H)_y \ \otimes \ {}_y (\cc \# H)_x \longrightarrow {}_z (\cc \# H)_x$$
is given by
$$(_z \varphi_y \otimes h) \circ(_y \psi_x \otimes h')=\sum {}_z\varphi_y \circ (h_1\ {}_y\psi_x) \otimes h_2h',$$
where the comultiplication $\Delta$ of $H$ is given by $\Delta (h)
= \sum h_1 \otimes h_2$ and $\circ$ denotes composition in $\cc$.
As before we have an immediate coherence result:

\begin{prop}
Let $\cc$ be a finite object $H$-module category  $\cc$. Then the
$k$-algebras $a(\cc)\# H$ and  $a(\cc \# H)$ are canonically
isomorphic.
\end{prop}

Let now $G$ be a group. A $G$-graded  $k$-category  $\cc$ (see for
instance \cite{cima}) is a $k$-category $\cc$ such that each
morphism space $_y\cc_x$ is the direct sum of sub-vector spaces
$_y\cc^s_x$, indexed by elements $s\in G$ such that ${_z\cc_y}^t \
{{}_y\cc_x}^s\subseteq {{}_y\cc_x}^{ts}$ for all $x,y \in \cc$ and
for all $s,t \in G$.

Notice that as in the algebra case, gradings of a $k$-category
$\cc$ by means of a group $G$ are in one-to-one correspondence
with $kG$-comodule category structures on $\cc$. Let now $G$ be a
finite group, $\cc$ be a $G$-graded $k$-category and consider the
function algebra $k^G=(kG)^*$ which is a Hopf algebra. The
category $\cc$ is a $k^G$-module category, hence according to our
previous definition we can consider $\cc \# k^G$.

We want to compare this category with another construction of a
$k$-category denoted $\cc \# G$ which can be performed for an
arbitrary group $G$, see \cite{cima} : the set of objects is
$\cc_0 \times G$ while the morphisms from $(x,s)$ to $(y,t)$ is
the vector space ${{}_y\cc_x}^{\left(t^{-1}s\right)}$. The
composition of morphisms is well-defined as an immediate
consequence of the definition of a graded category.

Notice that given a graded algebra $A$ considered as a single
object $G$-graded $k$- category, the preceding construction
provides a category with as many objects as elements of $G$, even
if $G$ is infinite. If $G$ is finite, the associated algebra is
known to be the usual smash product algebra $A\#k^G$, see
\cite{cima}.

We will recall below the definition of the module category of a
$k$-category in order to prove that in case of a finite group $G$
the module categories  over $\cc\#k^G$ and $\cc\#G$ are
equivalent.

First we introduce a general setting which is interesting by
itself.

\begin{defi}\label{contraction}
Let $\dd$ be a $k$-category equipped with a partition $E$ of the
set of objects $\dd_0$ by means of finite sets $\{E_i\}_{i\in I}$.
Then $\dd_E$ is a new $k$-category obtained by \emph{contraction}
along the partition, more precisely $I$ is the set of objects of
$\dd_E$ and morphisms are given by

$$_j(\dd_E)_i = \bigoplus_{y\in E_j\ \ x\in E_i} {}_y\dd_x. $$

Composition is given by matrix product combined with composition
of the original category. Notice that the identity map of an
object $i$ is given by $\sum_{z\in E_i} \hbox{ } _z1_z $, which
makes sense since $E_i$ is finite.
\end{defi}

\begin{ejem}  Let $A$ be an algebra and let $F$ be a complete finite family of orthogonal
idempotents in $A$ (we don't require that the idempotents are
primitive). Consider the category $\dd$ with set of objects $F$
and morphisms ${}_y\dd_{x} = yAx$. Then the contracted category
along the trivial partition with only one subset is a single
object category having endomorphism algebra $\bigoplus_{x,y\in
F}{} \hbox{}_{y}\dd_{x} = \bigoplus_{x,y \in F} yAx = A.$
\end{ejem}
We also observe that for a finite object $k$-category $\cc$, the
contracted category along the trivial partition is a single object
category with endomorphism algebra precisely $a(\cc)$. More
generally let $E$ be a partition of $\cc_0$, then the $k$-algebras
$a(\cc)$ and $a(\cc_E)$ are equal.

We will establish now a relation between $\dd$ and $\dd_E$ at the
representation theory level of these categories. In order to do so
we recall the definition of modules over a $k$-category.

\begin{defi} Let $\cc$ be a $k$-category. A left $\cc$-module $\M$ is a
collection of $k$-modules $\{_x\M\}_{x\in \cc_0}$ provided with a
left \emph{action} of the  $k$-modules of morphisms of $\cc$,
given by $k$-module maps $_y\cc_x \otimes_k {}_x\M \rightarrow
{}_y\M$, where the image of $_yf_x \otimes {}_xm$ is denoted
${}_yf_x\ _xm$, verifying the usual axioms:

\begin{itemize}
\item ${}_zf_y\ ({}_yg_x\ _xm) = \left({}_zf_y\ {}_yg_x\right)\
{}_xm$,

\item ${}_x1_x \ {}_xm = {}_xm.$
\end{itemize}
\end{defi}

In other words $\M$ is a covariant $k$-functor from $\cc$ to the
category of $k$-modules, the preceding explicit definition is
useful for some detailed constructions. We denote by $\cc-\Mod$
the category of left $\cc$-modules. In case of a $k$-algebra $A$
it is clear that $A$-modules considered as $k$-vector spaces
equipped with an action of $A$ coincide with $\mathbb{Z}$-modules
provided with an $A$-action. Analogously, $\cc$-modules as defined
above are the same structures than $\mathbb{Z}$-functors from
$\cc$ to the category of $\mathbb{Z}$-modules.

\begin{defi}
Two $k$-categories are said to be \emph{Morita equivalent} if
their left module categories are equivalent.
\end{defi}

\begin{prop}
Let $\dd$ be a $k$-category and let $E$ be a partition of the
objects of $\dd$ by means of finite sets. Then $\dd$ and the
contracted category $\dd_E$ are Morita equivalent.
\end{prop}

We notice that this result is an extension of the well known fact
that the category of modules over an algebra is isomorphic to the
category of functors over the category of projective left modules
provided by a direct sum decomposition of the free rank one left
module, obtained for instance through a complete system of
orthogonal idempotents of the algebra.

\proof Let $\M$ be a $\dd$-module and let $F\M$ be the following
$\dd_E$-module:
$$_iF\M= \bigoplus_{x\in E_i} {}_x\M\mbox{ for each $i\in I$},$$
the action of a morphism ${}_jf_i =(_yf_x)_{x\in E_i, y\in E_j}\in
{}_j(\dd_E)_i$ on ${}_im=(_xm)_{x\in E_i}\in{}_iF(\M)$ is obtained
as a matrix by a column product, namely:
$$_jf_i\ _im= (\sum_{x\in E_i}{}_yf_x\ _xm)_{y \in E_j}.$$
A $\dd-\Mod$ morphism  $\phi:\M\to \M'$ is a natural
transformation between both functors, {\it i.e.} a collection of
$k$-maps ${}_x\phi: {}_x\M \to {}_x\M'$, satisfying compatibility
conditions. We define $F\phi: F\M\to F\M'$ by:
$$_i(F\phi)=\bigoplus_{x\in E_i} {}_x\phi.$$

Conversely given a $\dd_E$-module $\N$, let  $G\N \in (\dd-\Mod)$
be the functor given by $_x(G\N)=e_x \left({}_i\N\right)$, where
$i$ is unique element in $I$ such that $x\in E_i$, and where $e_x$
is the idempotent $|E_i|\times|E_i|$ - matrix with one in the
$(x,x)$ entry and zero elsewhere.

The action of ${}_yf_x\in {}_y\dd_x$ on ${}_x(G\N)$ is obtained as
follows: let $i,j\in I$ be such that $x\in E_i$ and $y\in E_j$.
Let $\overline{{}_yf_x} \in {}_j(\dd_E)_i$ be the matrix with
${}_yf_x$ in the $(y,x)$ entry  and zero elsewhere. Then, for
$e_xn\in {}_x(G\N)$ we put $({}_yf_x)(e_xn)=
{}_j(\overline{{}_yf_x})_i\ {}_i(e_xn)\in
e_y\left({}_j\N\right)={}_y(G\N)$.

It is easy to verify that both compositions of $F$ and $G$ are the
corresponding identity functors.

\medskip

We will now apply the preceding result to the situation $\dd=
\cc\# G$ using the partition provided by the orbits of the free
$G$-action on the objects.

\begin{teo}
The $k$-categories $\cc \# G$ and $\cc \# k^G$ are Morita
equivalent.
\end{teo}
\proof We consider the contraction of $\cc\# G$ along the
partition provided by the orbits, namely for $x\in \cc_0$ we put
$E_x=\{(x,g)\ |\ g\in G\}$.  Observe that for all $x\in \cc_0$ the
set $E_x$ is finite since its cardinal is the order of the group
$G$. Moreover the set of objects  $\left((\cc\#G)_E\right)_0$ of
the contracted category is identified to $\cc_0$.

The morphisms from $x$ to $y$ in the contracted category are
$\bigoplus_{s,t \in G} {{}_y\cc_x}^{t^{-1}s}$. On the other hand
$$ _y(\cc \#k^G)_x = {}_y\cc_x \otimes k^G = \bigoplus_{v\in G} {{}_y\cc_x}^v \otimes k^G.$$
We assert that the contracted category $(\cc\#G)_E$ and $\cc\#
k^G$ are isomorphic. The sets of objects already coincide. We
define the functor $L$ on the morphisms as follows. Let
$_{(y,t)}f_{(x,s)}$ be an elementary matrix morphism of the
contracted category. We put
$$L\left(_{(y,t)}f_{(x,s)}\right) = f \otimes \delta_s \in  {{}_y\cc_x}^{t^{-1}s} \otimes k^G.$$
It is not difficult to check that $L$ is an isomorphism preserving
composition.

\begin{obser}
The categories $\cc\# G$ and $\cc\# k^G$ are not equivalent in
general as the following simple example already shows : let $A$ be
the group algebra $kC_2$ of the cyclic group of order two $C_2$
and let $\cc_A$ be the single object $C_2$-graded $k$-category
with $A$ as endomorphism algebra. The category $\cc\#C_2$ has two
objects that we denote $(*,1)$ and $(*,t)$, while $\cc\#k^{C_2}$
has only one object $*$. If $\cc\# G$ and $\cc\# k^G$ were
equivalent categories the
 algebras $End_{\cc\# C_2}((*,1))$ and $End_{\cc\#
k^{C_2}}(*))$ would be isomorphic. However the former is
isomorphic to $k$ while the latter is the four dimensional algebra
$End_{\cc\# k^{C_2}}(*)= (k \bigoplus kt)\otimes k^{C_2}$.
\end{obser}

\section{\sf $kG$-module categories}

Let $G$ be a group and let $\cc$ be a $kG$-module category. Using
the Hopf algebra structure of $kG$ and the preceding definitions
we are able to construct the smash category $\cc\# kG$. We have
already noticed that if $\cc$ is an object finite $k$-category
then the algebra $a(\cc \#kG)$ is the classical smash product
algebra $a(\cc)\# kG$.

According to \cite{cima} a $G$-$k$-category $\dd$ is a
$k$-category with an action of $G$ on the set of objects and, for
each $s\in G$, a $k$-linear map $s: {}_y\dd_x \to {}_{sy}\dd_{sx}$
such that $s(gf)=s(g)s(f)$ and $t(sf)=(ts)f$ for any composable
couple of morphisms $g, f$ and any elements $s, t$ of $G$. Such a
category is called a free $G$-$k$-category in case the action of
$G$ on the objects is a free action, namely the only group element
acting trivially on the category is the trivial element of $G$.

\begin{obser}
Notice that $kG$-module categories are $G$-$k$-categories
verifying that the action of $G$ on the set of objects is trivial.
\end{obser}

We need to associate a free $G$-$k$-category to a $kG$-module
category $\cc$, in order to perform the quotient category as
considered in \cite {cima}. For this purpose we consider
\emph{inflated categories} as follows.

\begin{defi}
Let $\cc$ be a $k$-category and let
$F=\left(F_x\right)_{x\in\cc_0}$ be a sequence of sets associated
to the  objects of $\cc$. The set of objects of the inflated
category $I_F\cc$ is
$$\left\{(x,i)\mid x\in\cc_0 \mbox{ and } i\in F_x\right\}$$
while ${}_{(y,j)}{(I_F\cc)}_{(x,i)}={}_y\cc_x$ with the obvious
composition provided by the composition of $\cc$. Alternatively,
consider $F$ as a map $\varphi$ from a set to $\cc_0$ such that
the fiber over each object $x$ is $F_x$. The set of objects of the
inflated category is the fiber product of $\cc_0$ with this set
over $\varphi$.

\end{defi}
\begin{obser}
Clearly an inflated category is equivalent to the original
category since all the objects with the same first coordinate are
isomorphic. Hence a choice of one object in each set $\{(x,i)\mid
i\in F_x\}$ provides a full sub-category of $I_{F}\cc$ which is
isomorphic to $\cc$.
\end{obser}
In case $\cc$ is a $kG$-module category we use the constant
sequence of sets provided by the underlying set of $G$. We obtain
a free action of $G$ on the objects of the inflated category
$I_{G}\cc$ by translation on the second coordinate. Moreover the
original action of $G$ on each morphism set of $\cc$ provides a
free $G$-$k$-category structure on the inflated category. More
precisely the $G$-action on the category $I_{G}\cc$ is obtained
through maps for each $u\in G$ as follows:
$$u: {}_{(y,t)}I_{G}\cc_{(x,s)} \to {}_{(y,ut)}I_{G}\cc_{(x,us)}$$
$$ u\left({}_{(y,t)}f_{(x,s)}\right)={}_{(y,ut)}(u\left({}_yf_x)\right){}_{(x,us)}.$$

As a next step we notice that the free $G$-$k$-category $I_{G}\cc$
has a skew category $(I_{G}\cc)[G]$ associated to it. In fact any
$G$-$k$-category has a related skew category  defined in
\cite{cima}. We recall that $(I_{G}\cc)[G]_0= (I_{G}\cc)_0 = \cc_0
\times G$. For $x,y \in \cc_0$ $t,s \in G$ we have
$${}_{(y,t)}(I_{G}\cc)[G]_{(x,s)} = \bigoplus_{u\in G} {}_{(y,ut)}(I_{G}\cc)_{(x,s)} =
\bigoplus_{u\in G} {}_y\cc_x = {}_y\cc_x \times G.$$

We are going to compare the categories $\cc \# kG$ and
$(I_{G}\cc)[G]$. In order to do so we consider the intermediate
quotient category $(I_{G}\cc)/G$ (see \cite[Definition
2.1]{cima}). We recall the definition of $\dd /G$, where $\dd$ is
a free $G$-$k$-category: the set of objects is the set of
$G$-orbits of $\dd_0$, while the $k$-module of morphisms in $\dd
/G$ from the orbit $\alpha$ to the orbit $\beta$ is
$$ _{\beta}(\dd /G)_{\alpha}= \left(\bigoplus_{b \in \beta, a\in \alpha} {}_b\dd_a\right)/G.$$
Recall that $X/G$ denotes the module of coinvariants of a
$kG$-module $X$, namely the quotient of $X$ by $(\ker\epsilon)X$
where $\epsilon:kG\to k$ the augmentation map. Composition is well
defined precisely because the action of $G$ is free on the
objects, more explicitly, for $g\in {}_d\dd_c$ and $f\in{}_b\dd_a$
where $b$ and $c$ are objects in the same $G$-orbit, let $s$ be
the unique element of $G$ such that $sb=c$. Then $[g][f]=[g\
(sf)]=[(s^{-1}g)\ f]$.

\begin{lem}\label{smashandquotient}
The $k$-categories  $\cc \# kG$ and $(I_{G}\cc)/G$ are isomorphic.
\end{lem}
\proof Clearly the set of objects can be  identified. Given a
morphism $({}_yf_x\otimes u)\in{}_y(\cc \# kG)_x$ we associate to
it the class $[f]$  of the morphism $f\in
{}_{(y,1)}(I_{G}\cc)_{(x,u)}$. Notice that in the smash category
we have
$$({}_zg_y\otimes v )({}_y f_x\otimes u) = {}_zg_y \ v ({}_y f_x)
\otimes vu$$ which has image $[{}_zg_y \ v ({}_y f_x)]$. The
composition in the quotient provides precisely $[g][f]=[g\ vf]$.
The inverse functor is also clear.

Since $(I_{G}\cc)/G$ and $(I_{G}\cc)[G]$ are equivalent (see
\cite{cima}), we obtain the following:

\begin{prop}
The categories $\cc \# kG$ and $(I_{G}\cc)[G]$ are equivalent.
\end{prop}

\section{\sf From Galois coverings to smash extensions and vice versa}

Our aim is to relate $kG$-smash extensions and Galois coverings
for a finite group $G$. Recall that it has been proved in
\cite{cima} that any Galois covering with group $G$ of a
$k$-category $\bb$ is obtained \textit{via} a $G$-grading of
$\bb$, we have that $\cc = \bb \#G$ is the corresponding Galois
covering of $\bb$. We have already noticed that for a finite group
$G$ a $G$-grading of a $k$-category $\bb$ is the same thing than a
$k^G$-module category structure on $\bb$.

However neither $\bb$ nor $\bb\# G$ have a natural $kG$-module
category structure which could provide a smash extension. We have
proven before that $\bb \# k^G$ is Morita equivalent to the
category $\bb \# G$. The advantage of $\bb\# k^G$ is that it has a
natural $kG$-module category structure provided by the left
$kG$-module structure of $k^G$ given by
$t\delta_s=\delta_{st^{-1}}$.

In this way we associate to the starting Galois covering $\bb\#G$
of $\bb$ the smash extension $(\bb\# k^G)\to (\bb\# k^G)\# kG$. In
\cite {rasa} the authors describe when a given Hopf-Galois
extension is of this type (in the case of algebras). We will prove
that the later is isomorphic to an \textit{ad-hoc} category
$M_{\mid G\mid}(\bb)$ which happens to be Morita equivalent to
$\bb$.

\begin{defi}
Let $\bb$ be a $k$-category and let $n$ be a sequence of positive
integers $\left(n_x\right)_{x\in\bb_0}$. The objects of the matrix
category $M_n(\bb)$ remain the same objects of $\bb$. The set of
morphisms from $x$ to $y$ is the vector space of $n_x$-columns and
$n_y$-rows rectangular matrices with entries in ${}_y\bb_x{}$.
Composition of morphisms is given by the matrix product combined
with the composition in $\bb$.
\end{defi}

\begin{obser}
In case the starting category $\bb$ is a single object category
provided by an algebra $B$, the matrix category  has one object
with endomorphism algebra precisely the usual algebra of matrices
$M_n(B)$.

Notice that the matrix category that we consider is not the
category $\mathrm{Mat}(\cc)$ defined by Mitchell in \cite{mi1}. In
fact $\mathrm{Mat}(\cc)$ corresponds to the additivisation of
$\cc$ (see the Appendix).
\end{obser}

We need the next result in order to have that the smash extension
associated to a Galois covering has categories Morita equivalent
to the original ones. In fact this result is also a categorical
generalization of Cohen Montgomery duality Theorem \cite{como}.

\begin{lem}
Let $\bb$ be a $G$-graded category. Then the categories $(\bb \#
k^G) \# kG$ and $M_n(\bb)$ are isomorphic.
\end{lem}
\proof Both sets of objects coincide. Given two objects $x$ and
$y$ we define two linear maps:
$$\phi: {}_y\bb_x \otimes k^G \otimes kG \to {{}_y}\!\left(M_n(\bb)\right)_x ,$$
$$\psi: {{}_y}\!\left(M_n(\bb)\right)_x \to {}_y\bb_x \otimes k^G \otimes kG.$$
Given an homogeneous  element $$(f\otimes \delta_g \otimes h) \in
{}_y\bb_x \otimes k^G \otimes kG,$$ where $f$ has degree $r$ and
$g,h \in G$ we put
$$\phi(f\otimes \delta_g \otimes h)= f \ \hbox{}_{rg}E_{gh},$$
where $_{rg}E_{gh}$ is the elementary matrix with $1$ in the
$(rg,gh)$-spot and $0$ elsewhere. It is straightforward to verify
that $\phi$ is well-behaved with respect to compositions.

We also define $\psi$ on elementary  morphisms as follows:
$$\psi(f\ {}_gE_h) = f\otimes \delta_{r^{-1}g} \otimes g^{-1}rh,$$
 where $r$ is the degree of $f$.

\bigskip

Next we have to prove that $M_n (\bb)$ is Morita equivalent to
$\bb$. In order to do so we develop some  Morita theory for
$k$-categories which is interesting by itself. When we restrict
the following theory to a particular object, it will coincide with
the classical theory, see for instance \cite[p.326]{we}. Moreover,
in case of a finite object set $k$-categories both Morita theories
coincide using the associated algebras that we have previously
described.

Let $\cc$ be a $k$-category. For simplicity for a given object $x$
we denote by $A_x$ the $k$-algebra $_x\cc_x$. For each $x$, let
$B_x$ be a $k$-algebra such that there is a  $(B_x, A_x)$-bimodule
$P_x$  and a $(A_x, B_x)$-bimodule $Q_x$ verifying that $P_x
\otimes_{A_x} Q_x \cong B_x$  as $B_x$-bimodules and $Q_x
\otimes_{B_x} P_x \cong A_x$  as $A_x$-bimodules. In other words
for each object we assume that we have a Morita context providing
that $A_x$ and $B_x$ are Morita equivalent. Note that it follows
from the assumptions that $P_x$ is projective and finitely
generated on both sides, see for instance \cite{we}.

Using the preceding data we modify the morphisms in order to
define a new $k$-category $\dd$ which will be Morita equivalent to
$\cc$. In particular the endomorphism algebra of each object $x$
will turn out to be $B_x$.

More precisely the set of objects of $\dd$ remains the set of
objects of $\cc$ while for morphisms we put
$$_y\dd_x = P_y \otimes_{A_y} {}_y\cc_x
\otimes_{A_x} Q_x.$$
 Notice that for $x = y$ we have $_x\dd_x
\cong B_x$. In order to define composition in $\dd$  we need to
provide a map
$$(P_z \otimes_{A_z} {}_z\cc_y \otimes_{A_y} Q_y)
\otimes_k (P_y \otimes_{A_y} {}_y\cc_y \otimes_{A_x} Q_x)
\longrightarrow P_z \otimes_{A_z} {}_z\cc_x \otimes_{A_x} Q_x,$$
For this purpose let $\varphi_x$ be a fixed $A_x$-bimodule
isomorphism from $Q_x\otimes_{B_x}P_x$ to  $A_x$ and consider
$\phi_x$ the composition the projection $Q_x\otimes_{k}P_x\to
Q_x\otimes_{B_x}P_x$ followed by $\varphi_x$. Then composition is
defined as follows
$$(p_z\otimes g \otimes q_y)(p_y\otimes f \otimes q_x)= p_z
\otimes g\left[\phi_y(q_y\otimes p_y)\right]f\otimes q_x.$$

This composition is associative since the use of the morphisms
$\phi$ do not interfere in case of composition of three maps.

\begin{prop}
Let $\cc$ and $\dd$ be $k$-categories as above. Then $\cc$ and
$\dd$ are Morita equivalent.
\end{prop}

\proof  For a $\cc$-module $\M$ we define the $\dd$-module $F\M$
as follows: $$_x(F\M) = P_x \otimes_{A_x} {}_x\M \mbox{, which is
already a left $B_x$-module.}$$ The left action $_y\dd_x \otimes
{}_x(F\M) \to {}_y(F\M)$ is obtained using the following morphism
induced by $\phi_x$
$$\left(P_y \otimes_{A_y} {}_y\cc_x \otimes_{A_x} Q_x \right) \otimes_k \left(P_x \otimes_{A_x}
{}_x\M\right) \longrightarrow  P_y \otimes_{A_y} {}_y\cc_x
\otimes_k {A_x}\otimes_k {}_x\M$$ and the actions of $A_x$ and of
$_y\cc_x$ on $_x\M$. We then obtain a map with target $_y(F\M)$.
This defines clearly a $\dd$-module structure.

Similarly we obtain a functor $G$ in the reverse direction which
is already an equivalent inverse for $F$.

\bigskip

We apply  now this Proposition to a $k$-category $\cc$ and the
category obtained from $\cc$ by replacing each endomorphism
algebra by matrix algebras over it. For each object $x$ in $\cc_0$
consider the $k$-algebra $B_x = M_n (A_x)$. The bimodule $_{M_n
(A_x)}(P_x)_{A_x}$ is the left ideal of $M_n (A_x)$ given by the
first column and zero elsewhere, while $_{A_x}(Q_x)_{M_n (A_x)}$
is given by the analogous right ideal provided by the first row.
Then the category $\dd$ defined above is precisely $M_n (\cc)$.

\begin{coro}
$\cc$ and $M_n (\cc)$ are Morita equivalent.
\end{coro}

\begin{obser}
  An analogous Morita equivalence still hold when the integer $n$
  is replaced by a sequence of positive integers
  $\left(n_x\right)_{x\in \cc_0}.$

\end{obser}

The applications of Morita theory for categories developed above
covers a larger spectra than the one considered in this paper. We
have produced several sorts of Morita equivalences for categories,
namely expansion, contraction and the Morita context for
categories described above. We will prove the next result in the
Appendix.

\begin{teo}\label{morita}
Let $\cc$ and $\dd$ be Morita equivalent $k$-categories. Up to
equivalence of categories, $\dd$ is obtained from $\cc$ by
contractions and expansions.

\end{teo}

\begin{ejem} Let $A$ be a $k$-algebra and $\cc_A$ the
corresponding single object category. It is well known that the
following $k$-category $M\cc_A$ is Morita equivalent to $\cc_A$:
objects are all the positive integers $[n]$ and the morphisms from
$[n]$ to $[m]$ are the matrices with $n$ columns, $m$ rows, and
with $A$ entries.

At each object [n] choose the system of $n$ idempotents provided
by the elementary matrices which are zero except in a diagonal
spot where the value is the unit of the algebra. The expansion
process through this choice provides a category with numerable set
of objects, morphisms are $A$ between any couple of objects, they
are all isomorphic, consequently this category is equivalent to
$\cc_A$. This way a Morita equivalence (up to equivalence) between
$\cc_A$ and $M\cc_A$ is obtained using the expansion process.

Conversely, in order to obtain $M\cc_A$ from $\cc_A$, first
inflate $\cc_A$ using the set of positive integers. Then consider
the partition by means of the finite sets having all the positive
integers cardinality, namely $\{1\}, \{2,3\}, \{4,5,6\}, \dots$.
Finally the contraction along this partition provides precisely
$M\cc_A$.
\end{ejem}

We provide now an alternative proof of the fact that a matrix
category is Morita equivalent to the original one. It provides
also evidence for Theorem \ref{morita} concerning the structure of
the Morita equivalence functors. First consider the inflated
category using the sequence of positive integers defining the
matrix category. We have shown before that this category is
equivalent to the original one. Secondly perform the contraction
of this inflated category along the finite sets partition provided
by couples having the same first coordinate. This category is the
matrix category. Since we know that a contracted category is
Morita equivalent to the original one, this provides a proof that
a matrix category is Morita equivalent to the the starting
category, avoiding the use of Morita contexts. The alternative
proof we have presented indicate how classical Morita equivalence
between algebras can be obtained by means of contractions,
expansions and equivalences of categories. More precisely Theorem
\ref{morita} states that classical Morita theory can be replaced
by those processes.

The results that we have obtained provide the following

\begin{teo} Let $\cc\longrightarrow\bb$ be a Galois covering of categories with finite group $G$.
The associated smash extension
$\bb\#k^G\longrightarrow(\bb\#k^G)\#kG$ verifies that  $\bb\#k^G$
is Morita equivalent to $\cc$ and $(\bb\#k^G)\#kG$ is Morita
equivalent to $\bb$.

\end{teo}

Finally notice that the proof of a converse for this result is a
direct consequence of the discussion we have made  in the previous
section:

\begin{teo}
Let $\cc\longrightarrow\bb$ be a smash extension with finite group
$G$. The corresponding Galois covering
$I_{G}\cc\longrightarrow\left(I_{G}\cc\right)/G$ verifies that
$I_{G}\cc$ is  equivalent to $\cc$ and that
$\left(I_{G}\cc\right)/G$ is  equivalent to $\bb$.

\end{teo}

\proof Indeed an inflated category is isomorphic to the original
one; moreover $\bb=\cc\#kG$ and by Lemma \ref{smashandquotient}
this category is isomorphic to $(I_{G}\cc)/G$.

\section{Appendix: Morita equivalence of categories over a field}

We have considered in this paper several procedures that we can
apply to a $k$-category. We briefly recall and relate them with
the karoubianisation (also called idempotent completion) and the
additivisation (or additive completion), see for instance the
appendix of \cite{ti}.

The \emph{inflation} procedure clearly provides an equivalent
category : given a set $F_x$ over each object $x$ of the
$k$-category $\cc$, the objects of the inflated category $I_F\cc$
are the couples $(x,i)$ with $i\in F_x$. Morphisms from $(x,i)$ to
$(y,j)$ remain the morphisms from $x$ to $y$. Consequently objects
with the same first coordinate are isomorphic in the inflated
category. Choosing one of them above each object of the original
category $\cc$ provides a full subcategory of the inflated one,
which is isomorphic to $\cc$.

The \emph{skeletonisation} procedure consists in choosing
precisely one object in each isomorphism set of objects and
considering the corresponding full subcategory. Clearly any
category is isomorphic to an inflation of its skeleton. Skeletons
of the same category are isomorphic, as well as skeletons of
equivalent categories.

Those remarks show that up to isomorphism of categories, any
equivalence of categories is the composition of a skeletonisation
and an inflation procedure.

Concerning Morita equivalence, we have used \emph{contraction} and
\emph{expansion}. In order to contract we need a partition of the
objects of the $k$-category $\cc$ by means of finite sets. The
sets of the partition become the objects of the contracted
category, and morphisms are provided by matrices of morphisms of
$\cc$. Conversely, in order to expand we choose a complete
 system of orthogonal idempotents for each endomorphism algebra at
each object of the $k$-category (the trivial choice is given by
just the identity morphism at each object).  The set of objects of
the expanded category  is the disjoint union of all those finite
sets of idempotents. Morphisms from $e$ to $f$ are $f{}_y\cc_xe$,
assuming $e$ is an idempotent at $x$ and $f$ is an idempotent at
$y$. Composition is given by the composition of $\cc$.

We assert that the karoubianisation and the additivisation (see
for instance \cite{basc,ti}) can be obtained through the previous
procedures.

Recall that the karoubianisation of $\cc$ replaces each object of
$\cc$ by \emph{all} the idempotents of its endomorphism algebra,
while the morphisms are defined as for the expansion process
above.

Consider now the partition of the objects of the karoubianisation
of $\cc$ given by an idempotent and its complement, namely the
sets $\{e, 1-e\}$ for each idempotent at each object of $\cc$. The
contraction along this partition provides a category equivalent to
$\cc$, since all the objects over a given object of $\cc$ are
isomorphic in the contraction of the karoubianisation. Concerning
the additivisation, notice first that two constructions are in
force which provide equivalent categories as follows.

The larger category is obtained from $\cc$ by considering all the
finite sequences of objects, and morphisms given through matrix
morphisms of $\cc$. Observe that two objects (i.e. two finite
sequences) which differ by a transposition are isomorphic in this
category, using the evident matrix morphism between them.

Consequently the objects of the smaller construction are the
objects of the previous one modulo permutation, namely the set of
objects are finite sets of objects of $\cc$ with positive integers
coefficients attached. In other words objects are maps from
$\cc_0$ to $\mathbb{N}$ with finite support. Morphisms are once
again matrix morphisms.

The observation above concerning finite sequences differing by a
transposition shows that the larger additivisation completion is
equivalent to the smaller one.

Finally the smaller additivisation of $\cc$ can be expanded:
 choose the canonical complete orthogonal idempotent system at
each object provided by the matrix endomorphism algebra. Of course
the expanded category have several evident isomorphic objects
which keeps trace of the original objects. A choice provides a
full subcategory equivalent to $\cc$.

It follows from this discussion that karoubianisation and
additivisation provide Morita equivalent categories to a given
category, using contraction and expansion processes, up to
isomorphism of categories.

We denote $\widehat{\cc}$ the completion of $\cc$, namely the
additivisation of the karoubianisation (or \emph{vice-versa} since
those procedures commute).
We notice  that two categories are
Morita equivalent if and only if their completions are Morita
equivalent.

Recall that a $k$-category is called \emph{amenable} if it has
finite coproducts and if idempotents split, see for instance
\cite{fr}. It is well known and easy to prove that the completion
$\widehat{\cc}$ is amenable.

We provide now a proof of Theorem \ref{morita}. We have shown that
the completion of a $k$-category is obtained (up to equivalence)
by expansions and contractions of the original one. Notice that
$\widehat\cc$ and $\widehat\dd$ are Morita equivalent amenable
categories. We recall now the proof that this  implies that the
categories $\widehat\cc$ and $\widehat\dd$ are already equivalent
(a result known as "Freyd's version of Morita equivalence", see
\cite[p.18]{mi1}): consider the full sub-category of representable
$\widehat\cc$-modules, namely modules of the form
${}_{-}\widehat\cc_x$. This category is isomorphic to the opposite
of the original one (this is well known and immediate to prove
using Yoneda's Lemma). Since $\widehat\cc$ is amenable,
representable $\widehat\cc$-modules are precisely the small (or
finitely generated) projective ones, see for instance \cite[p.
119]{fr}. Finally the small projective modules are easily seen to
be preserved by any equivalence of categories; consequently the
opposite categories of $\widehat\cc$ and $\widehat\dd$ are
equivalent, hence the categories themselves are also equivalent.

\footnotesize \noindent C.C.:
\\D\'epartement de Math\'ematiques,
 Universit\'e de Montpellier
2,  \\F--34095 Montpellier cedex 5, France. \\{\tt
Claude.Cibils@math.univ-montp2.fr}

\noindent A.S.:
\\Departemento de Matem\'atica,
 Facultad de Ciencias Exactas y Naturales,\\
 Universidad de Buenos Aires
\\Ciudad Universitaria, Pabell\'on 1\\
1428, Buenos Aires, Argentina. \\{\tt asolotar@dm.uba.ar}

\end{document}